\documentstyle[twoside,11pt,amssymb]{article}
\newtheorem{proposition}{Proposition}
\newtheorem{theorem}{Theorem}
\newtheorem{lemma}{Lemma}
\newtheorem{corollary}{Corollary}
\newtheorem{remark}{Remark}
\newtheorem{example}{Example}

\newtheorem{definition}{Definition}

\newcommand{\bth}{\begin{theorem}}
\newcommand{\bpr}{\begin{proposition}}
\newcommand{\epr}{\end{proposition}}
\newcommand{\bco}{\begin{corollary}}
\newcommand{\eco}{\end{corollary}}
\newcommand{\ble}{\begin{lemma}}
\newcommand{\ele}{\end{lemma}}
\newcommand{\bre}{\begin{remark}\rm}
\newcommand{\ere}{\end{remark}}
\newcommand{\bex}{\begin{example}\rm}
\newcommand{\eex}{\end{example}}
\newcommand{\bres}{{\bf Remarks}\begin{enumerate}}
\newcommand{\eres}{\end{enumerate}}
\newcommand{\bexs}{{\bf Examples}\begin{enumerate}\begin{enumerate}}
\newcommand{\eexs}{\end{enumerate}\end{enumerate}}
\newcommand{\bde}{\begin{definition}\rm}
\newcommand{\ede}{\end{definition}}
\newcommand{\bpf}{\begin{pf}}
\newcommand{\epf}{\end{pf}}
\newcommand{\bproof}{\begin{pf}}
\newcommand{\eproof}{\end{pf}}
\newcommand{\bpft}{\begin{prooft}\begin{itemize}}

\title{\huge
Some Remarks on \\ Symmetric Products of Curves}

\author{Sadok Kallel}

\begin{document}
\date{}
\maketitle



\parskip=0.5pc
\parindent=0pc
\font\sc=cmcsc10 at 11pt

\def\ra#1{\hbox to #1pc{\rightarrowfill}}
\def\fract#1#2{\raise4pt\hbox{$ #1 \atop #2 $}}
\def\decdnar#1{\phantom{\hbox{$\scriptstyle{#1}$}}
\left\downarrow\vbox{\vskip10pt\hbox{$\scriptstyle{#1}$}}\right.}
\def\decupar#1{\phantom{\hbox{$\scriptstyle{#1}$}}
\left\uparrow\vbox{\vskip10pt\hbox{$\scriptstyle{#1}$}}\right.}

\def\za{\vrule height6pt width4pt depth1pt}
\def\lrar{{\ra 2}}

\def\sp#1{SP^{#1}}
\def\spy{SP^{\infty}}

\def\tensor{\otimes}
\def\deg{\hbox{deg}}
\def\dim{\hbox{dim}}
\def\fil{\hbox{fil}}
\def\mod{\hbox{mod}}

\def\bbc{{\mathbb C}}
\def\bbp{{\mathbb P}}
\def\bba{{\mathbb A}}
\def\bbz{{\mathbb Z}}
\def\bbq{{\mathbb Q}}


\begin{abstract}
Symmetric products of curves are important spaces for both
geometers and topologists, and increasingly useful objects for physicists.
We summarize below some of their basic homotopy theoretic properties
and derive a handful of known and less well-known results about
them\footnote{We would like to thank
``x'' for correcting and improving an earlier version of this paper.}.
\end{abstract}

This paper combines in a slightly leisurely way both geometry and
topology to describe some useful properties of symmetric products of
algebraic curves. We record a straightforward derivation of Clifford's
theorem from a calculation of MacDonald, point to an equally simple
characterization of hyperelliptic curves and discuss the embeddability
(both in the continuous and holomorphic categories) of the unique
spherical generator in dimension two in the homology of these
spaces. A homotopy retract statement about the Abel-Jacobi map is also
proven.

\vskip 10pt \noindent{\bf\Large\S1 Cohomology Structure and
Clifford's Theorem}\vskip 10pt

Given a complex algebraic curve
$C$ and $n\geq 1$, the $n$-th symmetric product of $C$
is the quotient $C^{(n)}= C^n/\Sigma_n$, where $\Sigma_n$ is the
symmetric group acting on $C^n$ by permuting coordinates\footnote{In
the algebraic topology literature, this is
often written $\sp{n}(C)$.}. Elements in
$C^{(n)}$ are referred to as {\sl effective divisors} on $C$.  A
point $D\in C^{(n)}$ is said to have degree $n$ and we can write it
as a formal linear combination $\sum n_ix_i$ where $x_i\neq x_j\in C$ 
for $i\neq j$, and $n_i$ are positive integers with $\sum n_i=n$.

{\sc 1.1. MacDonald and Clifford's theorems}

It it assumed well-known that $C^{(n)}$ is a complex (smooth)
algebraic variety for all positive $n$.
The $n$-th Abel-Jacobi map is an algebraic map
$$\mu_n : C^{(n)}\lrar J(C)$$
where $J(C)$ is the ``Jacobian'' of $C$ (a complex torus of
dimension the genus of $C$). It is additive in the sense that
the following commutes
$$\matrix{C^{(r)}\times C^{(s)}&\lrar& C^{(r+s)}\cr
\decdnar{\mu_r\times\mu_s}&&\decdnar{\mu_{r+s}}\cr
J(C)\times J(C)&\lrar&J(C)\cr}$$
where the bottom map is addition in the abelian torus $J(C)$, and
the top map is concatenation of points (which is also an abelian pairing).
If $C$ is an elliptic curve for example,
then $J(C)\cong C$. The inverse preimages of $\mu_n$ are complex
projective spaces $\bbp^m$, where $m$ is related to the dimension of
some complete linear series on $C$(cf. \cite{cornalba}). The dimension
of the preimages $\mu_n^{-1}(x), x\in J(C)$
is an upper semi-continuous function of $x$.
When $n\geq 2g$, the dimension is constant for given $n$ and equals
$n-g$ (this is a direct consequence of Riemann Roch for curves). In
fact, Mattuck shows that $\mu_n$ is a projectivized analytic bundle
projection with fiber $\bbp^{n-g}$.
The fiber being Kahler and the coefficients being simple (as can be
checked), it follows by a theorem of Blanchard that this fibration
is homologically trivial and
$$H_*(C^{(n)};\bbz ) \cong H_*(\bbp^{n-g})\tensor H_*(S^1)^{\tensor 2g}
\ \ ,\ \ n\geq 2g$$
(here of course we use the fact that $J(C)\simeq (S^1)^{2g}$).
In fact we will show that $J(C)$ is a homotopy retract
of $C^{(n)}$ for $n\geq 2g$ (section 1.3).
The situation is less clear for $2\leq n < 2g$ since the fibers can
jump up in dimension. This jump is however well controlled
(see corollary \ref{clifford}).

The following main theorem of MacDonald describes the cohomology
of $C^{(n)}$ for all $n$. As it turns out, these spaces
have interesting cup products.

\bth\label{macdo} (I.G. Macdonald) Let $e_i^*\in H^1(C)$ be the one
dimensional generators, $1\leq i\leq 2g$ and $b^*\in H^2(C)$ the
orientation class. Then the cohomology ring of $C^{(n)}$ over the
integers ${\bf Z}$ is generated by the $e_i^*$ and $b^*$ subject to
the following relations:\hfill\break {(i)} The $e_i^*$'s anti-commute
with each other and commute with $b^*$;\hfill\break {(ii)} If $i_1,
\dots ,i_a, j_1,\dots ,j_b,k_1,\dots, k_c$ are distinct integers from
1 to $g$ inclusive, then $$ e_{2i_1-1}^*\cdots
e_{2i_a-1}^*e_{2j_1}^*\cdots e_{2j_b}^*(e_{2k_1-1}^*e_{2k_1}^*- b^*)
\cdots (e_{2k_c-1}^*e_{2k_c}^* - b^*) (b^*)^q = 0$$ provided that
$a+b+2c+q = n+1.$\hfill\break 
If $n < 2g$ all the relations above are
consequences of those for which $q=0,1$, and if $n > 2g-2$ all the
relations are consequences of the single relation
$$(b^*)^{n-2g+1}\prod_{i=1}^g(e_{2i-1}^*e_{2i}^*-b^*) = 0$$
\end{theorem}

In other words $H^*(C^{(n)})$ is the quotient of
$H^*(J;\bbz )[b^*] = E(e_1^*,\ldots, e_{2g}^*)\tensor \bbz [b^*]$
by the above relations, where $E$ stands for an exterior algebra.
In the first part of this note we describe the geometry behind
MacDonald's theorem and give some simple manipulations and
applications of the homology of symmetric products. For example a
rather pleasant application is the derivation, in
purely topological terms, of
a pivotal theorem in the geometry of algebraic curves
due to Clifford.

\bpr Let $C$ be a closed oriented topological surface of genus $g>0$,
and $\bbp^m\lrar C^{(n)}$ a map that is non-zero on
second homology groups. Id
\\ (i) $n< 2g$, then necessarily $m\leq
\left[{n\over 2}\right]$, \\ (ii) and if $n\geq 2g$, then $m\leq n-g$.
\epr

\noindent{\sc Proof:} Suppose $m = [\frac{n}2] + 1$ and $n < 2g$ and
suppose there is a continuous $h:\bbp^m\lrar C^{(n)}$
with $h^*(b^*) \not= 0$ (the induced map in cohomology).
Then necessarily $h^*(b^*)^m \not= 0$ by the ring
structure of projective space.  On the other hand, $h^*(e_i^*)=0$ and
hence
$$h^*(b^{*m}) = h^*( \prod_{k=1}^m (b^* - e_{ 2k-1}^* e_{2k}^*) ) = 0$$
using the relation of MacDonalds theorem \ref{macdo}. This contradicts the
choice of $m$. One uses a similar argument for (ii).
\hfill\za

\bco\label{clifford} (Clifford)
Suppose $0\leq n\leq 2g-1$, and pick $x\in J(C)$. Write
$\mu_n^{-1}(x) = \bbp^m$ for some $m$. Then necessarily
$m\leq {n\over 2}$.
\eco

\bre It is well-known that when $C$ is hyperelliptic, $C^{(2)}$ is the
blowup at a point of $W_2=\mu (C^{(2)})\subset J(C)$. The
exceptional fiber is $\bbp^1\subset C^{(2)}$ and hence for all $m$
one has an essential map $(\bbp^1)^{(m)}=\bbp^m\lrar C^{(2m)}$. 
The bound for
Clifford's theorem is clearly sharp.\ere

\bre
The class $\theta := \sum_1^ge_{2i-1}^*e_{2i}^*\in H^2(J;\bbz )$
is particularly
important since it is Poincar\'e dual to the class of the divisor
$W_{g-1}:=\mu [C^{(g-1)}]$. This (or a translate
of it) is known as the $\Theta$-divisor (\cite{cornalba}, chap I, \S4).
\ere

{\sc 1.2. Homology}

Alternatively, $C$ being a genus $g$ closed topological
surface, one can express the homology of $C^{(n)}$ as a ``truncated
Pontryagin ring" (as in \cite{k1} for instance). To be more explicit, write
$${\tilde H}_*(C;\bbz )\cong\cases{\bbz\{b\},& when $*=2$, with generator
$b$\cr \bbz^{2g}\{e_1,\ldots, e_{2g}\},& when $*=1$\cr}$$
The generators $e_i$ are spherical classes corresponding to the inclusion
of the one skeleton $\bigvee^{2g}_1S^1$ in $C$. It will be often convenient
to appeal to the cofibration sequence describing $C$
\begin{equation}\label{cofibration}
\bigvee^{2g}S^1\fract{i}{\lrar} C\fract{p}{\lrar} S^2
\end{equation}
We choose the $e_i$'s so that representative cycles of $e_{2i-1}$ and
$e_{2i}$ intersect transversally at a point ($e_i$ and $e_j$ are
disjoint otherwise). These generators have (canonical) duals $b^*$ and
$e_i^*$ ($b^*$ is the pull-back of the generator under the quotient
map $C\lrar S^2$).  This as we know implies that $H^*(C;\bbz )$ is
generated by $b^*, e_1^*,\ldots, e_{2g}^*$ with the relation $b^* =
e_{2i-1}^*e_{2i}^*$ (and graded commutativity).

The following, of which we give a short proof, is a well-known
property of symmetric products. Choose once and for all a basepoint
$x_0\in C$.

\ble \label{homologyembedding}
The basepoint inclusion $i_n: C^{(n)}\lrar C^{(n+1)}$, $D\mapsto
D+x_0$ induces a homology monomorphism.  \ele

\noindent{\sc Proof}: A class $\alpha\in H^i(C^{(n)})$ is determined
by a map $f: C^{(n)}\lrar K(\bbz, i)$ into an Eilenberg-MacLane
space. Composition gives a map $g: C\hookrightarrow C^{(n)}\lrar
K(\bbz, i)$. Since $K(\bbz, i)$ can be chosen to be an abelian
monoid\footnote{In fact $K(\bbz, i)$ is the infinite symmetric product
of the sphere $S^i$.}, we get a map $f\times g: C\times C^{(n)}\lrar
K(\bbz, i)$ and since this monoid is moreover commutative, we further
get a map $\beta: C^{(n+1)}\lrar K(\bbz, i)$.  This represents an
element in $H^i(C^{(n+1)})$ and by construction $i^*(\beta ) =
\alpha$. This shows that $i^*$ is surjective, and hence $i_*$ is
injective. \hfill\za

\bre
In fact a more general theorem of Steenrod (\cite{k2},
section 2.3) shows that the maps $i_*$ are split embeddings in homology.
This statement is valid for any connected simplicial complex $X$ and
an elementary illustration can be found in lemma \ref{circles}.
\ere

Note that the multiplication $C^{(r)}\times C^{(s)}\lrar C^{(r+s)}$
given by concatenation of points, has known effect in homology.
For instance the covering map $\pi : C^n\lrar C^{(n)}$ corresponds to the
iterated multiplication $(C^{(1)})^n\lrar C^{(n)}$, and under
$\pi$ the class $b^{\tensor n}$ maps to the orientation
class of $C^{(n)}$ (since this is an oriented manifold); that is
$$\pi_* : [C]^{\tensor n}\mapsto n![C^{(n)}]$$
where $n! = \deg \pi$.
It is convenient to write in this case $\gamma_n := [C^{(n)}]$. Write
$$\cdot : H_i(C^{(r)})\tensor H_j(C^{(s)})\lrar H_{i+j}
(C^{(r+s)})$$
the corresponding Pontryagin product. Under this product, we see that
$b\cdot b = 2\gamma_2$, and more generally we have

\bpr\label{homology}
\cite{k2}  $H_*(C^{(n)})$ has generators $\gamma_i =
[C^{(i)}]$, $e_i$, $1\leq i\leq 2g$, and all products of such
$e_{i_1}\cdot\ldots \cdot e_{i_k}\cdot \gamma_j$,
$i_1<i_2<\ldots <i_k$, of length at most $n$ (i.e. $k+j\leq n$).
The classes $\gamma$ verify $\gamma_i\cdot\gamma_j = \left({i+j\atop
i}\right)\gamma_{i+j}$.
When $n\geq 2g$, the exterior algebra
$E(e_1,\ldots, e_{2g})$ embeds into $H_*((C^{(n)})$.
 \epr

\bre
We do identify all throughout this paper $C^{(i)}$ with its image
in $C^{(n)}$ under the basepoint embeddings for $1\leq i\leq n$.
We also identify, as in proposition \ref{homology},
$H_*(C^{(i)})$ with its image in
$H_*(C^{(n)})$ so that for example the class $e_1\cdot
e_2\in H_2(C^{(3)})$ comes from the subspace $H_2(C^{(2)})$.
Technically it would have been better to say that
if $x\in H_*(C^{(i)})$ then $x\cdot\iota^{(n-i)}\in
H_*(C^{(n)})$ where $\iota$ is the generator in $H_0(C)$.
\ere

\bex
As an illustration, $H_*(C^{(2)})$ has generators $e_i$ (in dimension 1),
$b = [C]$ and all two fold products $e_i\cdot e_j, i< j$ (in dimension 2),
$e_i\cdot b$ in dimension 3 and then $b^2$ in dimension four.
The Euler characteristic is $\chi (S)=
(g-1)(2g-3)$ in agreement with the formula $\left( {{2g-2}\atop 2}\right)$
in \cite{macdo}.
\eex

There is a slick homotopy argument to see why proposition
\ref{homology} is true.  If we denote the direct limit of the
inclusions $i_n$'s (lemma \ref{homologyembedding})
by the infinite symmetric product $\spy (C)$ (\footnote{$\spy$
is a homotopy functor and we adopt here the notation of Dold and Thom.}),
then we have a monomorphism
$H_*(C^{(n)})\hookrightarrow H_*(\spy (C))$. Now $\spy (C)$ is a
very easy space to describe. This is the free abelian monoid on points
of $C$. We then have two maps:\\ $\bullet$ The map into the Jacobian
extends additively to $\mu : \spy (C)\lrar J(C)$.  \\ $\bullet$ There
is a map $p : \spy (C)\lrar\spy (S^2)$ induced from the quotient map
$C\lrar S^2$.\\
The crucial observation is now that

\ble \label{atinfinity}
$\mu\times p : \spy (C)\lrar \spy (S^2)\times J(C)\simeq
\bbp^{\infty} \times (S^1)^{2g}$ is a homotopy equivalence.  \ele

{\sc Proof}.
The composite $\bigvee S^1\hookrightarrow C\lrar J(C)$ induces an
isomorphism in $H_1$ by construction of the Abel-Jacobi map.
On the other hand, the map $C\lrar S^2$ in (\ref{cofibration})
sends orientation class
to orientation class. The main ingredient we need is the theorem
of Dold and Thom which asserts that $\pi_*(\spy (X))\cong
{\tilde H}_*(X;\bbz )$ for any finite type CW complex (cf.
\cite{k2}, chapter 2). Applying $\pi_*$ to both sides of
$\spy (C)\lrar \spy (S^2)\times J(C)$, we obtain the map
$H_*(C )\lrar H_*(S^2)\oplus H_*(\bigvee S^1)^{\tensor 2g}$
which by the first two observations is an isomorphism. It follows
that $\mu\times p$ is a weak equivalence, and hence a homotopy
equivalence.
\hfill\za

{\sc Proof of proposition \ref{homology}}.  The equivalence $\spy
(C)\fract{\simeq}{\lrar} J(C)\times\bbp^{\infty}$ constructed above is
an equivalence of $H$-spaces since the map is multiplicative.  It
follows that $H_*(\spy (C),\bbz )$ (as a ring under the pontryagin
product $\cdot$) is an exterior algebra on $2g$ generators $e_i$ of
dimension 1 (the homology of $J(C)$), tensored with a divided power
algebra $\Gamma [b] =
\bbz\{b=\gamma_0, \gamma_1,\ldots,\}$ on a 2-dimensional generator $b
= [C]$ (which is the homology of $\bbp^{\infty}$). By construction
$H_*(C^{(n)})$ includes in this ring as those products of
lengths at most $n$. \hfill\za

{\sc 1.3 The Map into the Jacobian}

The map $C^{(n)}\lrar J(C)$ is injective in cohomology for all
$n\geq 2g$. This is because the composite
$$\bigvee^{2g} S^1\hookrightarrow C\lrar C^{(n)}\fract{\mu}{\lrar} J(C)$$
is an isomorphism on $H_1(-;\bbz )$ for all $n\geq 1$, and hence for
$n\geq 2g$, the image of $\mu^*$ is precisely the exterior algebra on
the $e_i^*$ and is an isomorphism onto its image.

\bre
To see this injectivity statement in a more conceptual level, one can
quote a general fact pointed out by Gottlieb \cite{gottlieb}. Let $f:
M\lrar N$ be a smooth surjective map between oriented closed
manifolds, and suppose $F=f^{-1}(y)$ (for some regular value $y\in N$
and hence $F$ is a closed manifold) is not a boundary. Then $f^*:
H^*(N;\bbz_2)\lrar H^*(M;\bbz_2)$ is injective. Note that in our case
when $f=\mu$, $F$ is a projective space and hence cannot be a
boundary.  Moreover since the spaces in questions have no torsion,
Goettlieb's mod-$2$ statement is true integrally. Note that the Abel-Jacobi
map is no longer smooth for $n < 2g$.\ere

\bpr \label{retract}
$J(C)$ is a homotopy retract of $C^{(n)}$ when $n\geq 2g$.
\epr

By that we mean there are maps $J(C)\lrar C^{(n)}\fract{\mu}{\lrar}
J(C)$ of which composite is homotopic to the identity of $J(C)$.  This
was obviously the case when $n=+\infty$ according to lemma
\ref{atinfinity}, so this is represents a refinement to finite $n$.
The statement is in fact a painless corollary of the following useful
but not difficult fact (this is part of the main theorem of \cite{ong}
with a longer argument of proof based on the theory of hyperplane
arrangements).

\ble\label{circles} There is a homotopy equivalence
$\left(\bigvee^kS^1\right)^{(n)} \simeq (S^1)^k$ for $n\geq k\geq 1$.
\ele

{\sc Proof}. When $k=1$, this is easy since $(S^1)^{(n)}$ can be
identified with $(\bbc^*)^{(n)}$, and the map
$$ (\bbc^*)^{(n)}\lrar {\mathcal H},\ \ 
\sum z_i\mapsto\prod_{1\leq i\leq n}(z-z_i)
$$ 
is a homeomorphism, where on the left we have identified the space of
monic polynomials of degree $n$ with $\bbc^n$, and those polynomials
with roots avoiding the origin with $\mathcal H$;
the complement of a hyperplane in
$\bbc^n$. Since ${\mathcal H} \cong\bbc^{n-1}\times\bbc^*\simeq S^1$, the
claim follows in this case.
It is amusing to see for instance that $(S^1)^{(2)}$ is
precisely the mobius band (cf. \cite{k2}, chapter 2).  We next assume
$n\geq k >1$.

The map inducing the homotopy equivalence is simply given by
the embedding into the infinite symmetric product
$$\phi : (\bigvee^kS^1)^{(n)}\lrar \ \spy (\bigvee^kS^1)
= \prod^k\spy (S^1) \simeq (S^1)^k$$
Here we use the fact that $\spy (S^1)\simeq S^1$, or more generally
that $\spy (S^n)$ is a model for $K(\bbz,n)$. 
Note that the composite 
$\bigvee^kS^1\lrar (\bigvee^kS^1)^{(n)}\fract{\phi}{\lrar} 
(S^1)^k$ is an injection in homology according to lemma
\ref{homologyembedding}, and an isomorphism on $H_1$ between two
copies of $\bbz^{k}$.
Since symmetrization abelianizes fundamental groups as soon as 
$n\geq 2$, and since $H_1$ is $\pi_1$ made abelian, it follows that
$\phi$ induces an isomorphism in $\pi_1$ as well.

The homology of $(S^1)^k$ is obviously exterior on one-dimensional
generators $E(e_1,\ldots, e_k )$. Let $e_i \in H_1((\bigvee S^1)^{(n)})$
be the generators in dimension one coming from the various factors.
Using lemma \ref{atinfinity}, and the fact that
$\phi$ is a multiplicative map, we find that
$H_*((\bigvee^k S^1)^{(n)})$ (once identified with its image in the exterior
algebra $E(e_1,\ldots, e_k)$) is generated by all monomials of the form
$$e_{i_1}\tensor\cdots\tensor e_{i_r}\ ,\ \ \sum i_j \leq n$$
But when $n\geq k$, these account for \emph{all} generators in 
$E(e_1,\ldots, e_k)$ meaning that $\phi$ is actually an isomorphism
in homology. But then the spaces are simple ($(S^1)^k$ being a group
even) and so by Whitehead's theorem, the map $\phi$ is a homotopy
equivalence.

Note the one thing we used above is that $H_*(\spy (\bigvee^k
S^1))\cong E(e_1,\ldots, e_k )$ not only as abelian groups but also
and more importantly as Pontryagin rings.  \hfill\za

Proposition \ref{retract} now follows by setting $k=2g\leq n$ and 
observing that the composite\footnote{
Given $f : X\rightarrow Y$, we define its symmetrization
$f^{(n)} : X^{(n)}\lrar Y^{(n)}$ in the obvious way 
$f^{(n)} ( \sum n_ix_i) := \sum n_if(x_i)$.}
$$J(C)\simeq (\bigvee^{2g}S^1)^{(n)}
\fract{i^{(n)}}{\ra 4}C^{(n)}\fract{\mu}{\lrar} J(C)$$
is precisely the homotopy equivalence $\phi$ indicated above, with
$i : \bigvee^{2g}S^1\hookrightarrow C$ is as in (\ref{cofibration}).

{\sc 1.4 Intersection Theory on Symmetric Products}

Since $M=C^{(n)}$ is a $2n$-dimensional (oriented) manifold, we have
an intersection pairing $\bullet$ which relates to cup product via
Poincar\'e duality $pd$ as follows
$$\matrix{H_i(M)\tensor H_j(M)&\fract{\bullet}{\lrar}&H_{i+j-2n}(M)\cr
\decdnar{pd}&&\decdnar{pd}\cr
H^{2n-i}(M)\tensor H^{2n-j}(M)&\fract{\cup}{\lrar}&H^{4n-i-j}(M).
\cr}$$
When $i=j=n$ and upon identifying $H_0(M)\cong\bbz$,
$\bullet$ becomes the intersection from in middle homology.
If $x$ is a (co)cycle, we denote by $x^{\perp}$ its Poincar\'e
dual. For example $C^{(n-i)}\hookrightarrow C^{(n)}$ is a subvariety
of codimension $2i$ and we claim that
\begin{equation}\label{dualofdivisor}
[C^{(n-i)}]^{\perp} = (b^*)^i
\end{equation}
To see this we rely on the following principle (which is
discussed in \cite{gottlieb} for example) :
Suppose $f:M\lrar N$ is a (smooth) map between two closed oriented manifolds.
Then if a cycle and a cocycle are related by
Poincar\'e duality in $N$, $f^{-1}$(cycle) and $f^*$(cocycle) are
related by Poincar\'e duality in $M$. Alternatively if $X$ is a cycle on $N$
such that $f^{-1}(X)$ is defined as a cycle, and $Y$ is a cycle on $M$ such
that $Y\bullet f^{-1}(X)$ is defined, then in homology
$$f(Y\bullet f^{-1}(X)) = f(Y)\bullet X$$
(the right side beeing automatically defined). Let's apply this principle
to the projection $\pi : C^n\lrar C^{(n)}$. The inverse image of
the cycle $[C^{(n)}]$ is $\sum [C]^i$ where the sum is over all
possible ways to include $[C]^{n-i}:=b^{\times n-i}$ into $b^{\times n}$.
But the Poincar\'e dual
to $b^{\times n-i}$ is $(b^*)^{\times i}$ inserted in the remaining spots.
For example the dual to $b\tensor 1$ in $H_*(C^2)$ is $1\tensor b^*$.
And so we see that $(\pi^{-1}[C^{(n-i)}])^{\perp}$ coincides
with $\pi^*(b^*)^i$, and by the above principle
$[\sp{n-1}]^{\perp} = (b^*)^i$.

For the sake of the next argument,
write $e_i$ and $e_{i+g}$ the symplectic duals in $H_1(C)$ (instead of
$e_{2i-1}$ and $e_{2i}$). We can then show in $H_*(C^{(n)})$ that
\begin{eqnarray*}
e_i^{\perp} &=& (b^*)^{n-1}e_{i+g}^*
\ \ \ \ , \ \ \   (b\cdot e_i)^{\perp} = (b^*)^{n-2}e_{i+g}^*
\end{eqnarray*}
To see the first equality for example, notice that in $C^n$,
$1\times\cdots 1\times e_{i}\times 1\cdots 1$ is Poincar\'e dual to
$b^*\times\cdots b^*\times e_{i+g}^*\times\cdots\times b^*$.
And so
$$[\pi^{-1}(e_i)]^{\perp} =\left[\sum 1\tensor e_i\tensor 1\right]^{\perp}
= \sum_r (b^*)^r\tensor e_{i+g}^*\tensor (b^*)^{n-r-1}
= \pi^*\left((b^*)^{n-1}e_{i+g}^*\right)$$
With a little more care one can show that

\ble $(e_{i}\cdot e_j)^{\perp} =
\cases{(b^*)^{n-2}e_{i+g}^*e_{j+g}^*& $j\neq i+g\ \mod\ (2g)$\cr
(b^*)^{n-1} + (b^*)^{n-2}e_{i+g}^*e_{j+g}^*,& $j = i+g$ mod $(2g)$\cr}$
\ele

The difference in expression is of course due to the fact that
$e_i^*e_j^*=0$ if $i\not\equiv j+g\ \mod (2g)$ in $H^2(C)=\bbz$,
and that $e_i^*e_{i+g}^*=1$.
The following corollary is useful for later and can of course be deduced
form \cite{macdo}

\bco\label{signature} The signature of $S:=C^{(2)}$
is $\sigma (S) = 1-g$.
\eco

{\sc Proof}. We write down the intersection matrix! A basis for
$H_2(C^{(2)})$ is given by $b$ and then all possible classes
$e_i\cdot e_j$ with $i< j$. We have that $b^{\perp} = b^*$,
 $(e_i\cdot e_{i+g})^{\perp} = b^* - e_i^*e_{i+g}^*$ and that
$(e_i\cdot e_j)^{\perp} = e_{i+g}^*e_{j+g}^* = -e_i^*e_j^*$
if $1\leq i < j \neq i+g \leq 2g$.
Upon identifying $H_0=H^4=\bbz$ (and $(b^*)^2=1$), we see that $b\bullet b =1$.
For $j= i+g $, $s = r+g$, we see that
$$(e_i\cdot e_j)\bullet (e_r\cdot e_s)=
(b^*+e_{i+g}^*e_{j+g}^*)(b^*+e_{r+g}^*e_{s+g}^*) =
(b^*-e_{i}^*e_{j}^*)(b^*-e_{r}^*e_{s}^*) =
-1$$
according to theorem \ref{macdo}. This shows that the self-intersection of
$(e_i\cdot e_j)$ is trivial if $i\not\equiv j\mod (2g)$ and is
$-1$ otherwise. Since there are $g$ generators of the form
$e_ie_{i+g}$ with self-intersection $-1$, and a single generator $b$ with
self-intersection $+1$, the claim follows.
\hfill\za

\bre (Chern Classes). The total Chern class of $C^{(n)}$ is computed
in \cite{macdo}, p.332, from which one deduces in particular that
\begin{equation}\label{chernclass}
c_1 = (n-g+1)b^* - \theta\ \ ,\ \
c_2 = (n-g+1)(n-g).{(b^*)^2\over 2} - (n-g)b^*\theta + {\theta^2\over 2}
\end{equation}
By using Hirzebruch
signature formula $\sigma (S) = {1\over 3}<c_1^2-2c_2, [S]>$, the fact
that $x\theta  =g (b^*)^2$ in
$C^{(2)}$ and that $<(b^*)^2,[C^{(2)}]>=1$,
we recover corollary \ref{signature} immediately.
\ere


\vskip 10pt \noindent{\bf\Large\S2 The Spherical Class}\vskip 10pt

Generally and for any curve $C$, one easily constructs an essential map
$\alpha : S^2\lrar C^{(2)}$ as follows.  One can see directly (and for
any $X$) that
$\pi_1(X^{(n)})\cong H_1(X;\bbz )$ for $n\geq 2$
(cf. \cite{k2}, proposition 2.4 for a very short
proof). Now $C$ as a CW complex is a
two disk $D^2$ attached to the bouquet $\bigvee^{2g}S^1$
via a product of commutators in the free group
$\pi_1(\bigvee^{2g} S^1)$.
The attaching map (on $\partial D^2=S^1$)
becomes null-homotopic when mapping into $C^{(2)}$ since
$\pi_1(C^{(2)})$ is abelian.  It follows that {up to homotopy} the
basepoint embedding $C\lrar C^{(2)}$ factors through the cofiber $S^2$
and we have a map
$$C\lrar S^2\fract{\alpha}{\ra 2} C^{(2)}$$

\ble\label{sphericalclass}
$u:=\alpha_*([S^2]) = b - \ell$, where $\ell=
\sum_{i\leq g} e_{2i-1}\cdot e_{2i} \in H_2(C^{(2)}).$ \ele

\noindent {\sc Proof:} \cite{k1}
$h_*([S^2])$ is spherical and hence primitive.
We will be done if we show that the only primitive classes in
dimension 2 are ( multiples) of $b - \ell$.  The classes of dimension
two in $H_*(C^{(n)})$ are $b$ and the products $e_i\cdot e_j$ and so we
need determine the diagonal on each. For $b=[C]$ this is fairly
direct. Assume first that $C=T$ is a torus; i.e. $T\simeq S^1\times
S^1$, with one-dimensional generators $e_1$ and $e_2$. These classes
are primitive with $\Delta_*(e_i)=e_i\tensor 1+1\tensor e_i$ and hence
\begin{eqnarray*}
\Delta_*([T])&=&\Delta_*(e_1\tensor e_2)=(e_1\tensor 1+1\tensor e_1)
(e_2\tensor 1+1\tensor e_2)\\ &=&[T]\tensor 1 + e_1\tensor e_2 -
e_2\tensor e_1 + 1\tensor [T]
\end{eqnarray*}
When $C$ is of genus $g\geq 1$, we have a map $C\lrar T_1\vee\ldots
\vee T_g$ which is an isomorphism on $H_1$ and maps $[C]$ to $\oplus
[T_i]$. Combining these facts yields
$$\Delta_*([C])=[C]\tensor 1 +\sum_{i\leq g} e_{2i-1}\tensor
e_{2i}-\sum e_{2i} \tensor e_{2i-1}+1\tensor [C]$$ On the other hand
by tracing through the commutative diagram
$$\matrix{C\times C&\fract{\Delta\times\Delta}{\ra 3}~ C^2\times C^2~
\fract{1\times T\times 1}{\ra 3}&C^2\times C^2\cr
\decdnar{\cdot}&&\decdnar{\cdot \times \cdot }\cr C^{(2)}&\fract{\Delta}{\ra
8}&C^{(2)}\times C^{(2)}\cr} $$
(where $\cdot$ refers to the symmetric product pairing
and $T$ is the
appropriate shuffle map), we readily determine that
$$\Delta_*(e_i\cdot e_j) = e_i\cdot e_j\tensor 1 + e_i\tensor
e_j-e_j\tensor e_i + 1\tensor e_i\cdot e_j$$ Let $\ell=\sum_{1\leq
i\leq g} e_{2i-1}\cdot e_{2i}$. By the formulae above
$\bar\Delta_*\left([C]-\ell\right)= 0$ (where $\bar\Delta$ is the
reduced diagonal), hence implying that $[C]-\ell$ is primitive, and by
inspection the unique such (up to multiple). To show this multiple is +1,
and since the embedding $C\hookrightarrow\sp{2}(C)$ is homotopic
to $g : C\lrar S^2\lrar \sp{2}(C)$ we must have
$g^*(u^*) = b^*$. But then if $u=k(b-\ell )$,
$g^*(u^*)= kg^*(b^*-\ell^*) = kg^*(b^*)=kb^*$. The claim follows.
\hfill\za

\ble We have $(e_{2i-1}\cdot e_{2i})^* = e_{2i-1}^*e_{2i}^* - b^*$, with
the product on the left being the symmetric product pairing and on the
right the cup product. \ele

\noindent{\sc Proof:} We wish to determine the hom-dual of
$e_{2i-1}\cdot e_{2i}$ in $H^2(C^{(n)})$. By an earlier remark,
$\pi^*: H^*(C^{(n)})\lrar H^*(C)^{\tensor n}$ is injective,
where $\pi : C^n\lrar C^{(n)}$ is the covering projection. The
class $e_{2i-1}\cdot e_{2i}$ is the image of classes of the form $\pm
1\tensor\cdots 1\tensor e_r\tensor 1\cdots 1\tensor e_s\tensor
1\cdots\tensor 1$ with $\{r,s\}=\{2i-1, 2i\}$ (it is + if $r<s$ and
$-$ if $r>s$). It follows that
$$\pi^*(e_{2i-1}\cdot e_{2i})^* = \sum \pm 1\tensor\cdots 1\tensor
e_r^*\tensor 1\cdots 1\tensor e_s^*\tensor 1\cdots\tensor 1,\ \
\{r,s\}=\{2i-1, 2i\}$$
The sum is over all possible distinct spots $e_r$ and $e_s$ can assume.  On the
other hand, $\pi^*(e_i^*) = \sum 1\tensor\cdots 1\tensor e_i^*\tensor
1\cdots\tensor 1$ and hence
\begin{eqnarray*}
\pi^*(e_{2i-1}^*e_{2i}^*) &= &\pi^*(e_{2i-1}^*)\pi^*(e_{2i}^*)\\
&=&\pi^*(e_{2i-1}\cdot e_{2i})^* + \sum 1\tensor\cdots\tensor
e_{2i-1}^*e_{2i}^* \tensor 1\cdots\tensor 1\\
&=&\pi^*(e_{2i-1}\cdot e_{2i})^* + \sum 1\tensor\cdots\tensor b^* \tensor
1\cdots\tensor 1\\ &=&\pi^*(e_{2i-1}\cdot e_{2i})^* + \pi^*(b^*)
\end{eqnarray*}
Since $\pi^*$ is injective we get the desired equality. \hfill\za

\bre More generally we can show that
$(e_{2k_i-1}^*e_{2k_i}^*- b^*)
\cdots (e_{2k_c-1}^*e_{2k_c}^* - b^*)$ is dual to the $2c$-fold product
$e_{2k_i}\cdot e_{2k_i}\cdot\cdots e_{2k_c-1}\cdot e_{2k_c}$.
Evidently if the length of this class which is $2c$ exceeds $n$, then
by lemma \ref{homology}, the class is trivial in $H_{2c}C^{(n)}$ and
this explains MacDonald's relation
$(e_{2k_i-1}^*e_{2k_i}^*- b^*)\cdots (e_{2k_c-1}^*e_{2k_c}^* - b^*)=0$
for $2c>n$.
\ere

\ble Suppose $\alpha$ is a differential embedding, and denote by
$u^{\perp}$ the Poincar\'e dual of $u:=\alpha_*[S^2]\subset
H_2(C^{(2)})$. Then $u^{\perp} = (1-g)b^* + \theta$.
\ele

\noindent{\sc Proof}:
According to lemma \ref{sphericalclass},
$u^{\perp} = b^{\perp} - \left(\sum e_{2i-1}\cdot e_{2i}\right)^{\perp}$
The claim follows form the calculations in \S1.3.
Note that in light of lemma 3, the hom-dual and the Poincar\'e
dual of $u$ coincide.  \hfill\za

\bco\label{eulerclass}
If $\alpha$ is an embedding, then the euler class of the normal
bundle to $\alpha (S^2)$ is $1-g$ where $g$ is the genus of $C$.  \eco

\noindent{\sc Proof}: One need compute the self-intersection of
$\alpha (S^2)$ in $C^{(2)}$. This is of course given by the
evaluation $\langle (u^{\perp})^2, [C^{(2)}]\rangle$. We again
identify $H^4=\bbz$ so that $(b^*)^2 = e_{2i-1}^*e_{2i}^*e_{2j-1}^*e_{2j}^*=1$
for $i\neq j$. From this we get that $\theta^2 = g(g-1)$ and hence
\begin{eqnarray*}
(u^{\perp})^2 = \left[(1-g)b^* + \theta\right]^2
&=& (1-g)^2(b^*)^2 + 2(1-g)b^*\theta + \theta^2\\
&=& (1-g)^2 + 2(1-g)g + g(g-1) = (1-g)
\end{eqnarray*}
as claimed.\hfill\za

Note that when $\alpha : \bbp^1\lrar C^{(2)}$
is a holomorphic embedding, the normal bundle to the image (a Riemann
sphere) is a complex line bundle which is according to the above (and
as is well-known) isomorphic to ${\cal O}(1-g)$. In fact holomorphic
embeddings of $\alpha$ gives a characterization of hyperelliptic
curves.

\bpr\label{sphereinhyperelliptic}
The generator $u$ in $H_2(C^{(2)})$ can be represented by
a holomorphically embedded sphere (i.e. a rational curve)
if and only if $C$ is hyperelliptic.
\epr

\noindent{\sc Proof}: The ``only if'' part is straightforward since if
$\bbp^1$ is holomorphically embedded in $C^{(2)}$, then the
composite $\bbp^1\hookrightarrow C^{(2)}\lrar J(C)$ is trivial
necessarily and $C^{(2)}$ carries a $g_2^1$ \cite{cornalba}; that is
$C$ hyperelliptic.
Conversely, if $C$ hyperelliptic with a degree two branched covering
$\psi : C\lrar\bbp^1$, then the transfer map $\beta : \bbp^1\lrar
C^{(2)}$, sending $x$ into the unordered pair of preimages in
$\psi^{-1}(x)$, is an embedding of $\bbp^1$ into $C^{(2)}$.  The
class $\beta_*[S^2]$ is then a multiple of $u\in H_2(C^{(2)})$.
This multiple $k$ is determined by the intersection multiplicity of
$\beta (\bbp^1 )$ with the image under $\pi$ of $p\times C$ in
$C^{(2)}$ and this is of course a single point. So if we write
$\beta_*[\bbp^1 ]^{\perp}=ku^{\perp}$, the preceding discussion
shows that $ku^{\perp}b^* = 1$. But $ku^{\perp}b^* = k[(1-g)b^*+\theta]b^*
= k[(1-g) + g]=k$. The claim follows.
\hfill\za

\vskip 5pt {\sc 2.1. Rational Curves and the Symmetric Square} 

The following is very classical (algebraically) but we give it our own
spin (topologically).  For a hyperelliptic curve
$C\fract{:2}{\lrar}\bbp^1$, the transfer map $\bbp^1\lrar C^{(2)}$ is
a holomorphic embedding, and is homologous to the class $u$.  One
might then ask whether there are any other rational curves in
$C^{(2)}$ of degree $k>1$? Here $C$ is not necessarily hyperelliptic.
We can address this question using very classical geometry.  First
observe that an algebraic curve $X$ (of genus $g$) in $S:=C^{(2)}$ is
(by definition) a divisor on an algebraic surface and hence gives rise
to a line bundle $E$ on $S$ (cf. \cite{shavarevich}, chapter VI, \S1)
with the property that $E_{|X} = N_{X}$ (i.e. the restriction of $E$
to $X$ is the normal bundle of $X$ in $S$). An application of
Riemann-Roch in this situation shows the following \emph{adjunction}
formula: let $g_X$ be the genus of $X$, then
$$g_{X} = {(X+K)\bullet X\over 2} + 1$$
where $K$ is the canonical class of $S$, and where $\bullet$ the
intersection pairing (which is well defined within a divisor class,
\cite{shavarevich}). For example, when $X=\bbp^2$ then $K = -3L$,
where $L$ the divisor class containing all straight lines in $\bbp^2$.

\bco Identify $H_2(C^{(2)};\bbz )$ with $\bbz$. If the class $k$ 
is represented by a rational curve, then necessarily $k=1$ and
$C$ is hyperelliptic.\eco

{\sc Proof}.
The Poincar\'e dual of the canonical divisor for $C^{(n)}$ has been
computed by MacDonald (\cite{macdo}, p: 334) and it coincides with
\begin{equation}\label{firstchern}
K^{\perp} = -c_1(C^{(n)}) = (g-n-1)b^* + \sum e_{2i-1}^*e_{2i}^*
= (g-n-1)b^* + \theta
\end{equation}
Set $n=2$. Using the fact that $b^*=b^{\perp}$ and the formula
for $u^{\perp}$, this can be rewritten as
$K^{\perp} = (g-3)b^{\perp} + u^{\perp} - (1-g)b^{\perp}
= (2g-4)b^{\perp}+u^{\perp}$, and since the Poincar\'e duality is an
isomorphism, we find that
\begin{equation}\label{canonical}
K = (2g-4)b + u\ \ , \ \ K\in H_2(C^{(2)})
\end{equation}
Assume $X$ is a holomorphic sphere in $C^{(2)}$ with
$[X] = ku\in H_2(C^{(2)})$. Then by the adjuntion
formula $0 = (K+X)\bullet X + 2$.
Using that $b\bullet u=1$, $u\bullet u=1-g$, we write
$$-2 = (K+X)\bullet X = [(2g-4)b+ku ]\bullet ku =
k(2g-4) + (1+k)k(1-g)$$
This can be rewritten in the form $k(g-1)[1-k] = 2[k-1]$ from which
we deduce that the only possibility is when
$k=1$. According to proposition \ref{sphereinhyperelliptic}, this
implies the claim.\hfill\za

\bre An interesting observation (Mattuck) is that the canonical class
$K$ in the case $n=g$ is the union of all the special fibers of
$C^{(g)}\lrar J(C)$.  When $g=2$, the curve is automatically
hyperelliptic and $\mu$ has a single exceptional fiber $\bbp^1$.  In
this case $K=u$ the spherical class, in agreement with the
calculation (\ref{canonical}).  \ere


\vskip 10pt \noindent{\bf\Large\S3 Embedded Spheres in $C^{(2)}$}\vskip 10pt

The question now is whether any other
multiple of $u\in H_2(C^{(2)};\bbz )\cong\bbz$, can be realized by a
\emph{differentiably} embedded $2$-sphere in $C^{(2)}$. This we settle
partially as follows.

\bpr \label{embedding} Let $\beta :S^2\hookrightarrow C^{(2)}$ be an
embedding, with $\beta_*[S^2]  = ku$
in $H_2(C^{(2)})=\bbz$. If $g=0$ then $k=\pm 1, \pm
2$, and if $k>1$ odd and $g>0$ even, then necessarily $k\equiv\pm 1\
\mod\ (8)$. \epr

\bre For $g=0$ the answer is known \cite{lawson} and the following
cute result is attributed to Tristam. Let $u\in H_2(\bbp^2)$ be the
generator. Then $ku$ can be represented by an embedded $2$-sphere if
and only if $|k|<3$. The sufficiency is clear since the conic
$\bbp^1\hookrightarrow\bbp^2$, $[z,w]\mapsto [z^2:wz:w^2]$ is an
embedding of degree $2$. There are no other holomorphic embeddings of
higher degree since by the degree-genus formula an algebraic curve of
genus $g$ in the projective plane must satisfy $g = {1\over
2}(k-1)(k-2)$ and so it is rational if and only if $k=1,2$ as pointed
out. In the differentiable case, the proof is much less evident and
can be seen as a consequence of the theorem of Mrowka-Kronheimer
(previously the Thom conjecture) which asserts that any homology class
of degree $k$ in $\bbp^2$ can be smoothly realized by a curve of genus
at least ${1\over 2}(k-1)(k-2)$.\ere

A condition for embedding a $2$-sphere into a four manifold (not necessarily
simply connected) was first given by Kervaire and Milnor. It goes as follows.
An integral homology class $\beta\in H_2(M,\bbz )$ is
called ``characteristic'' if it is dual to the Stiefel-Whitney class
$w_2(M)$; i.e. if mod-$2$ reduction followed by Poincar\'e duality
takes $\beta$ to $w_2(M)$. Equivalently if the intersections verify
\begin{equation}
\beta\bullet \tau\equiv \tau\bullet \tau\ \mod (2)
\end{equation}
for any $\tau\in H_2(M)$. Note that any homology
class $\tau\in H_2(M;\bbz_2 )$ can be represented by an embedded
(real) surface (Thom). With integral coefficients the situation is
entirely different. The following is the main criterion of
Kervaire-Milnor

\bth\cite{kervaire}
Let $\beta\in H_2(M;\bbz )$ be dual to $w_2(M)$ where $M$ is a
closed, connected, oriented differentiable (real) $4$-manifold. Then
$\beta$ can be represented by a smoothly embedded 2-sphere only if
$\beta\bullet \beta\equiv\sigma (M)\ \mod\  16$.
\end{theorem}

Here $\sigma (M)$ is the signature of the intersection form.
Proposition \ref{embedding} is now a direct corollary.

\noindent{\sc Proof of Proposition \ref{embedding}}: Set $M=S$ the
second fold symmetric product of $C$. The class
$w_2(S)$ is the reduction mod-$2$ of the first chern class (for
unitary bundles). We've seen (\ref{firstchern}) that
$c_1(S) = (3-g)b^* - \sum e_{2i-1}^*e_{2i}^*$ and hence
reducing mod-$2$, one gets
$w_2(S) = b^* + \sum (b^* + e_{2i-1}^*e_{2i}^*)$ which is
precisely the mod-2 Poincar\'e dual of $u =
\alpha_*([S^2])$ (according to lemma 4). So $u$ is
characteristic and so is any odd multiple.
There now remains to compute the signature of $S$ and this is given
in lemma \ref{signature}.
Given then $S^2\lrar C^{(2)}=M$ which corresponds in homology to the
class $\beta = ku$, we have that $\beta\bullet \beta = k^2(1-g)$ (this is
the euler class of the normal bundle of $\beta$ if it's an embedding
according to corollary \ref{eulerclass}), and hence by the
Kervaire-Milnor congruence if $\beta$ is differentiably embedded (and
$k$ odd), then necessarily
$$k^2(1-g)\equiv 1-g\ \mod (16)$$
When $k=1$ the
sphere always embeds as we previously argued. When $k>1$ and $g=0$,
there is no odd $k$ such that the congruence is satisfied (as we know
already from Thom's conjecture). For even positive genus, the
congruence gives the restriction that $16\ |\ k^2-1$. Since $k$ odd,
this implies that $k\equiv\pm 1\ \mod\ (8)$.
\hfill\za

\vskip 5pt {\sc 3.1 The Second Homotopy Group}

We don't know if embeddings occur in those cases listed in proposition
\ref{embedding}, but we suspect the answer is no. Since the homotopy
class of $\alpha$ is in $\pi_2(C^{(2)})$ so we have better determine
what this group is.

\ble Suppose $C$ is of genus $g$, and $n > g$, then $\pi_2(C^{(n)})
\cong \bbz$.  \ele

{\sc Proof}.  When $n\geq 2g-1$, a theorem of Mattuck asserts that
$C^{(n)}$ fibers over $J(C)$ (the Jacobian) with fiber $\bbp^{n-g}$.
Since $\pi_3(J(C))=\pi_2(J(C))=0$ (being a complex torus), and since
$\pi_2\bbp^i = \bbz$ for $i\geq 1$, the claim follows in this
range. For $g\leq n\leq 2g-2$, the Abel-Jacobi map $\mu_n$ is a
quasifibration\footnote{A quasifibration $E\lrar B$ with ``fiber" $F$ has
a slightly weaker homotopy lifting property than fibrations but still enjoys
one of its main preperties : the long exact sequence of homotopy groups holds
for quasifibrations.} up to dimension $2(n-g)+1$ according to (\cite{guest},
theorem 4.1). So for $n>g$, it is a quasifibering up to dimension $3$
and so (by definition) there is a short exact sequence
$$\pi_3(J(C))=0\lrar \pi_2(\bbp^{n-g})\lrar \pi_2(C^{(2)})\lrar
\pi_2(J(C))=0$$
and as in the previous case this gives $\bbz$ for an answer.\hfill\za

Now and perhaps surprisingly $\pi_2(C^{(n)})$ for $n\leq g$ is not
even necessarily finitely generated (as pointed out to us by ``x'').
The map $C^{(n)}\lrar J(C)$ in this range has a description of a
blowup over various loci in the Jacobian of varying codimension. In
particular there is always an exceptional $\bbp^1\hookrightarrow
C^{(n)}$ and this generates a cyclic group in $\pi_2$ (the reason
being that the composite
$$\pi_2(\bbp^1)\fract{\cong}{\lrar} H_2(\bbp^1)\hookrightarrow H
_2(C^{(n)})$$ is essential), see \cite{bertram}. 
There might be however an interesting
action of $\pi_1$ which, unlike the case $n>g$, is not trivial.

\ble Suppose $C$ is a curve of genus $g=2$. Then there is a Laurent
polynomial description
$$\pi_2(C^{(2)}) =
\bbz [t_i, t_i^{-1}], \ 1\leq i\leq 4$$
\ele

{\sc Proof}.
The curve $C$ being of genus two, it is necessarily hyperelliptic
and so it is obtained from its Jacobian (a complex $2$-torus) by
blowing up a single point. The universal cover of the torus is
$\bbc^2$ and $J(C)$ is the quotient by a lattice $\mathcal L$
whose vertices are in one-to-one correspondance with $\bbz^4$.
It follows that the universal cover $\tilde X$ of
$C^{(2)}$ is a blowup
of $\bbc^2$ at each and everyone of these vertices. Each exceptional
fiber of this blowup being a copy of $\bbp^1$, it contributes a
generator to $H_2(\tilde X;\bbz )$. But then
$$ H_2(\tilde X) = \pi_2(\tilde X) = \pi_2(C^{(2)})$$
and so $\pi_2(C^{(2)})$ is infinitely generated indeed.
Choose a fundamental domain for the lattice made up of four
essential edges mapping to the fundamental group generators of
$J(C)$, which we call $t_1,\ldots, t_4$ and which also correspond to
the generators of $\pi_1(C^{(2)})=\bbz^4$ since the Abel-Jacobi map is
an isomorphism on $\pi_1$. Moving from generator to generator of
$H_2(\tilde X)$ through the lattice corresponds in $C^{(2)}$ to
multiplication by a word in the $t_i$'s or their inverses.
\hfill\za


\vskip 10pt

\small
\addcontentsline{toc}{section}{Bibliography}
\bibliography{biblio}
\bibliographystyle{plain}

\vskip 10pt
Sadok Kallel\\
Universit\'e Lille I, France\\
{\sc Email:} sadok.kallel@math.univ-lille1.fr

\end{document}